\documentclass[graybox]{svmult}

\usepackage{mathptmx}       
\usepackage{helvet}         
\usepackage{courier}        
\usepackage{type1cm}        
\usepackage{makeidx}         
\usepackage{graphicx}        
\usepackage{multicol}        
\usepackage[bottom]{footmisc}
\usepackage{amsmath,amssymb,amsfonts,color}
\usepackage{algorithmic}
\usepackage{algorithm,booktabs}
\usepackage[mathscr]{eucal}
\usepackage{subfigure}
\usepackage{xspace}
\usepackage{array}

\makeindex             
                       
\DeclareMathOperator*{\argmin}{arg\,min}

\begin{document}

\title*{Data-Driven Modeling and Control of Complex Dynamical Systems Arising in Renal Anemia Therapy}
\titlerunning{Data-Driven Control of Dynamical Systems in Renal Anemia Therapy}
\author{S. Casper, D.H. Fuertinger, P. Kotanko, L. Mechelli, J. Rohleff and S. Volkwein}
\institute{Sabrina Casper, Doris H. Fuertinger \at Fresenius Medical Care Deutschland GmbH, Bad Homburg, Germany. \\ \email{\{Sabrina.Rogg,Doris.Fuertinger\}@fmc-ag.com}	
\and Peter Kotanko \at
Renal Research Institute New York, New York, USA.
\\ \email{Peter.Kotanko@rriny.com}
\and Luca Mechelli, Jan Rohleff, Stefan Volkwein \at Department of Mathematics and Statistics, University of Konstanz, Konstanz, Germany.  \\ \email{\{luca.mechelli,jan.rohleff,stefan.volkwein\}@uni-konstanz.de}}
%
%
\maketitle

\abstract{This project is based on a mathematical model of erythropoiesis for anemia \cite{Mechelli_MS23_Fue12,Mechelli_MS23_FKTLK13}, which consists of five hyperbolic population equations describing the production of red blood cells under treatment with epoetin-alfa (EPO). Extended dynamic mode decomposition (EDMD) is utilized to approximate the non-linear dynamical systems by linear ones. This allows for efficient and reliable strategies based on a combination of EDMD and model predictive control (MPC), which produces results comparable with the one obtained in \cite{Mechelli_MS23_RFVKK19} for the original model.}
\vspace{-2em}
\section{Introduction}
\vspace{-1em}
\label{Mechelli_MS23_sec:1}
Almost all hemodialysis patients suffer from chronic anemia, due to the reduced functionality of the kidneys and the resulting low production of erythropoietin, a kidney-derived hormone that increases red blood cell output by the bone marrow. Therefore, physicians use erythropoietin stimulating agents, such as epoetin-alfa (EPO), to partially correct the anemia. The challenge in designing efficient therapies is due to the patients' differences in long-term response to EPO. In \cite{Mechelli_MS23_FKTLK13}, the authors introduce a mathematical model for predicting such a response. As in \cite{Mechelli_MS23_RFVKK19}, our aim is to design a feedback control strategy, based on Model Predictive Control (MPC) \cite{Mechelli_MS23_GP17}, to optimize the injections of EPO doses in order to reach a target hemoglobin level. In contrast to \cite{Mechelli_MS23_RFVKK19}, we do not imply the EPO model from \cite{Mechelli_MS23_FKTLK13} during the optimization, but we utilize it to generate a data-driven approximation of such a model through Extended Dynamic Mode Decomposition (EDMD) \cite{Mechelli_MS23_WKR15}. We then use such an EDMD-based surrogate model during the MPC. We refer to \cite{Mechelli_MS23_KM2018} where the authors introduced the idea to combine MPC and EDMD. We organize the work as follows: in Section~\ref{Mechelli_MS23_sec:2} we introduce the model and how we combine EDMD and MPC. In Section~\ref{Mechelli_MS23_sec:3} we show the reliability of the technique comparing our algorithm with the one in \cite{Mechelli_MS23_RFVKK19}.
\vspace{-2em}
\section{The EDMD-MPC based algorithm}
\vspace{-1em}
\label{Mechelli_MS23_sec:2}
Let us consider a finite horizon time interval $[0,T]$ with $T\gg 0$ and a finite number of injections $n_u\in\mathbb{N}$. We consider discrete injections of doses $u_i$ at predifined injection times $t_i$ for $i=1,\ldots,n_u$. Given a set of EPO doses $\textbf{u}=(u_i)_{i=1}^{n_u}\in U_\mathsf{ad}= \left\{\textbf{u}\in \mathbb{R}^{n_u}| 0\leq u_i\leq u_\text{max}, 1\leq i\leq n_u\right\}$, following \cite{Mechelli_MS23_RFVKK19}, one can compute the corresponding EPO concentration in the blood, which we indicate with $E(t;\textbf{u})$. This nonlinear function is continuously differentiable in $[0,T]$ and twice continuously differentiable in $U_\mathsf{ad}$. For further details see \cite{Mechelli_MS23_FKTLK13,Mechelli_MS23_RFVKK19}. The EPO model is composed of five coupled advection-reaction partial differential equations (PDEs) of the form
\begin{equation}
\label{Mechelli_MS23_EPOmodel}
\begin{aligned}
y_t(t,x) &= \kappa(x,E(t;\textbf{u}))y(t,x)-v(E(t;\textbf{u}))y_x(t,x) && \text{in } Q, \\
y(t,\underline{x}) &= g(t;E(t;\textbf{u})) && \text{in } (0,T), \\
y(0,x) &= y_0(x) && \text{in } \Omega
\end{aligned}
\end{equation}
with the spatial domain $\Omega=(\underline{x},\overline{x})\subset \mathbb{R}$, the space-time cylinder $Q=(0,T)\times\Omega$ and the initial condition $y_0$. The solution $y(t, x)$ to \eqref{Mechelli_MS23_EPOmodel} denotes the cell density of the respective cell population with maturity $x$ at time $t$. The coupling among the five equations is hidden in the boundary value $g(t;E(t;\textbf{u}))$. For the sake of brevity, we omit further explanations and we refer to \cite{Mechelli_MS23_FKTLK13,Mechelli_MS23_RFVKK19}. Note that in \eqref{Mechelli_MS23_EPOmodel}, the control $\textbf{u}$ (EPO doses) enters in the equations through the nonlinear EPO concentration, on which the advection and reaction coefficients depend. Such a complicated relationship between the states (cell densities) $y$ and the control $\textbf{u}$ leads to a non-convex optimization problem. Our aim is then to introduce a linear surrogate model based on EDMD, such that the resulting control-to-state map is linear. Doing so, we have the advantage of obtaining a convex optimization problem, which admits a unique minimizer and can be solved faster. In contrast to the standard Dynamic Mode Decomposition (DMD) \cite{Mechelli_MS23_KBBP16}, EDMD utilizes snapshots of the dynamic, controls and observables of a dynamical system to extract a discrete surrogate model \cite{Mechelli_MS23_WKR15}. Further we discretize in space with Legendre polynomials (as \cite{Mechelli_MS23_RFVKK19}) and in time with a constant time step $\Delta t$. In what follows $n$ is the number of Legendre polynomials and $m$ is the number of time steps. Let $\psi: \mathbb{R}^n \to \mathbb{R}^{N_\phi}$ be a vector of lifting functions (or observables) $\psi_i: \mathbb{R}^n\to \mathbb{R}$, i.e. $\psi(x) =\left(\psi_1(x),\ldots,\psi_{N_\phi}(x)\right)\in \mathbb{R}^{N_\phi}$. Let $Y_0= [y_0|...|y_{m-1}]\in\mathbb{R}^{n\times m}$ and $Y_1= [y_1|...|y_m]\in\mathbb{R}^{n\times m}$ and $U = [u_0|\ldots|u_{m-1}]\in\mathbb R^{1\times m}$ be given snapshots data matrices. Next, we need to define the matrices
\[
Y_{0,\text{lift}} = \left[
\begin{array}{ccc}
\psi_1(y_0)&\ldots&\psi_1(y_{m-1})\\
\vdots&&\vdots\\
\psi_{N_\phi}(y_0)&\ldots&\psi_{N_\phi}(y_{m-1})
\end{array}
\right], \hspace{2pt} 
Y_{1,\text{lift}} = \left[
\begin{array}{ccc}
\psi_1(y_1)&\ldots&\psi_1(y_m)\\
\vdots&&\vdots\\
\psi_{N_\phi}(y_1)&\ldots&\psi_{N_\phi}(y_m)
\end{array}
\right]\in\mathbb{R}^{N_\phi\times m}
\]
and identify the matrices $A\in\mathbb{R}^{N_\phi\times N_\phi}$, $B\in\mathbb{R}^{N_\phi\times n_u}$ and $C\in\mathbb{R}^{n\times N_\phi}$ such that
\begin{equation}
\label{Mechelli_MS23_argmin}
[A,B]= \argmin_{\tilde A,\tilde B}{\|Y_{1,\text{lift}}- \tilde AY_{0,\text{lift}}- \tilde BU\|}_F, \quad C=\argmin_{\tilde C}{\|Y_0-\tilde CY_{0,\text{lift}}\|}_F. 
\end{equation}
To solve \eqref{Mechelli_MS23_argmin} numerically, we have to perform two singular value decomposition; cf. \cite{Mechelli_MS23_KM2018}. We get then a discrete linear dynamical system in the observable space
\begin{equation}
\label{Mechelli_MS23_EDMDsyst}
\begin{aligned}
z_{k+1} &= Az_k + Bu_k\text{ for }k\ge0,\quad z_0 =\left(\psi_1(y_0),
\ldots,\psi_{N_\phi}(y_0)\right), \\
\hat{y}_{k+1} &= C z_{k+1}\text{ for }k\ge0,
\end{aligned}
\end{equation}
where $\hat{y}$ is the EDMD approximation of the $y$ solution to \eqref{Mechelli_MS23_EPOmodel}. Note that reconstructing and solving the EDMD system \eqref{Mechelli_MS23_EDMDsyst} is generally cheaper than solving the system of coupled hyperbolic PDEs \eqref{Mechelli_MS23_EPOmodel}. Thus, \eqref{Mechelli_MS23_EDMDsyst} is the surrogate model we will use during the MPC algorithm. More specifically, our goal is to apply EDMD to the fifth (and last) equation of the EPO model \eqref{Mechelli_MS23_EPOmodel} for two different choices of control snapshots, i.e. the EDMD matrix $U$ will contain:
\begin{enumerate}
	\item[(EDMD-C)\hspace{-5em}] \hspace{5em} The continuous EPO concentration $E(t;\textbf{u})$ at each time step;
	\item[(EDMD-D)\hspace{-5em}] \hspace{5em} The discrete EPO doses $\textbf{u}$ at the injection times and $0$ for the rest.
\end{enumerate}
The cost functional for the model in \cite{Mechelli_MS23_RFVKK19} discretized by Legendre polynomials is
\begin{align*}
J_N(\tilde y,\textbf{u}) & = \frac{1}{2} \sum_{j=1}^{n_u} \gamma_j u_j^2 + \frac{\sigma_\omega}{2}\int_{0}^{T} \left( \omega_5^{1/2}(\tilde y_5)_0(t)-P^\mathsf{d}\right)^2 \,\mathrm{d} t \\ & \quad + \frac{\sigma_f}{2}\left(\omega_5^{1/2}(\tilde y_5)_0(T)-P^\mathsf{d}\right)^2,
\end{align*}
where $(\cdot)_0$ means first-component in space and all the other parameters can be found in \cite{Mechelli_MS23_RFVKK19}. Moreover, $\tilde y = \mathcal{S}_N(\textbf{u})$ is the solution of (24b) in \cite{Mechelli_MS23_RFVKK19}.
Note that $J_N$ depends only on the solution of the fifth (and last) equation of \eqref{Mechelli_MS23_EPOmodel}. Therefore, to build our EDMD model, we use the solution $\tilde y_5$ as snapshots of the dynamic. We then get matrices $A$, $B$ and $C$ for \eqref{Mechelli_MS23_EDMDsyst} according to the chosen snapshots. Applying a trapezoidal rule to $J_N$ and considering the EDMD approximation $\hat{y}$, one obtains a cost functional $J_{m}(z,\textbf{u})= J_{N,\text{disc}}(Cz,\textbf{u})=J_{N,\text{disc}}(\hat{y},\textbf{u})$ and the optimal control problem is given as
\begin{align}
\label{Mechelli_MS23_PEDMD}
\min J_m(z,\textbf{u})\quad\text{s.t.}\quad z=\{z_k\}_{k=0}^m\subset\mathbb R^{N_\phi}\text{ satisfy \eqref{Mechelli_MS23_EDMDsyst}} \text{ for } \textbf{u}\in U_\mathsf{ad}.
\end{align}
Note that \eqref{Mechelli_MS23_PEDMD} is a convex linear-quadratic optimal control problem and thus admits a unique minimizer \cite{Mechelli_MS23_HPUU09}. If the EDMD approximation is sufficiently good, such a minimizer is also close to a local one of the optimal control problem in \cite{Mechelli_MS23_RFVKK19}. Since the horizon $T$ is generally large, the corresponding time discretization $t_0=0,\ldots,t_m= T$ contains many points. To avoid costly computations and compute a feedback control we introduce an MPC \cite{Mechelli_MS23_GP17} framework. This technique consists in fixing a prediction horizon $M$ and computing the solution of a first optimal control problem in $[0,M\Delta t]$, then storing the optimal control for the first time step, applying it to \eqref{Mechelli_MS23_EPOmodel} and repeating the procedure for the horizon $[t_1,t_1+M\Delta t]$ and so on. This approach has also the advantage of obtaining a feedback control which reacts to the solution of the EPO model. Since the EDMD is just a local approximation of the dynamics, we need to define a strategy in order to update the EDMD model during the MPC. We simply measure the difference between the EDMD solution and its EPO model counter part for the first time step. Note that this does not require additional computations with respect to the described procedure. We resume our algorithm in Algorithm~\ref{Mechelli_MS23_Alg1}.
\setlength{\textfloatsep}{1em}
\begin{algorithm}[t]
	\begin{algorithmic}[1]
		\STATE \textbf{Data:} Initial control $\textbf{u}^{(0)}\in U_\mathsf{ad}$, initial condition $y_\text{MPC}^{(0)}$, EDMD update tolerance $\tau_\text{upd}$, MPC prediction horizon $M$, EDMD update steps $M_\text{EDMD}$.
		\STATE Initialize the EDMD model at $\textbf{u}^{(0)}$ with snapshots from \eqref{Mechelli_MS23_EPOmodel};
		\FOR {$i=0,1, \ldots $} {
			\STATE Solve \eqref{Mechelli_MS23_PEDMD} in $[t_i, t_i+M\Delta t]$ with projected BFGS and initial guess $y_\text{MPC}^{(i)}$ to get an optimal $\bar{\textbf{u}}^{(i)}$;
			\STATE Store $u_\text{MPC}(t_i) = \bar{\textbf{u}}^{(i)}(t_i)$ and compute the new initial guess $y_\text{MPC}^{(i+1)}$ from \eqref{Mechelli_MS23_EPOmodel} using $u_\text{MPC}(t_i)$;
			\IF {$\|y_\text{MPC}^{(i+1)}-y_\text{EDMD}^{(i+1)}\|>\tau_\text{upd}$} 
			\STATE Update the EDMD model using snapshots of \eqref{Mechelli_MS23_EPOmodel} for $M_\text{EDMD}$ steps with the control $\bar{\textbf{u}}^{(i)}$;
			\ENDIF
		}
		\ENDFOR 
	\end{algorithmic}
	\caption{(EDMD-MPC algorithm) \label{Mechelli_MS23_Alg1}}
\end{algorithm}
\vspace{-2em}
\section{Numerical experiments}
\vspace{-1em}
\label{Mechelli_MS23_sec:3}
In this section, we compare the nonlinear MPC from \cite{Mechelli_MS23_RFVKK19} with the linear MPC based on the EDMD approach. For brevity, we report only the results on test conducted on Patient 2 and 3 of \cite{Mechelli_MS23_RFVKK19}. Let us mention, that we performed the numerical experiments also for the other patients in \cite{Mechelli_MS23_RFVKK19} and we got approximation errors similar to the one presented below. For both EDMD-MPC models, we choose an update tolerance $\tau_\text{upd} = 0.01$, $M_\text{EDMD}= 30$ steps, the MPC prediction horizon $M=14$ days and
\[
	\psi_1(y) = \omega_5^{1/2} (y)_0, \quad \psi_{i,j}(y) = L^j\left(\frac{(y)_i}{\sum_{i=1}^n (y)_i}\right), \quad i=0,\ldots,n, j=0,\ldots N_L-1,
\]
where $(y)_i$ is the $i$-th component of $y$ and $L^j$ is the $j$-th Legendre polynomial, $N_L = 2$ for EDMD-C and $N_L=6$ for EDMD-D.
All the other parameters are chosen as in \cite[Tables~4-8]{Mechelli_MS23_RFVKK19}. First of all, we consider a total time period of three weeks, i.e. $T=21$ days. We have an injection at day 1, 3 and 5 of each week.
\begin{figure}
	\vspace{-0.5em}
	\sidecaption
	\centering
	\includegraphics[width=0.49\textwidth]%
	{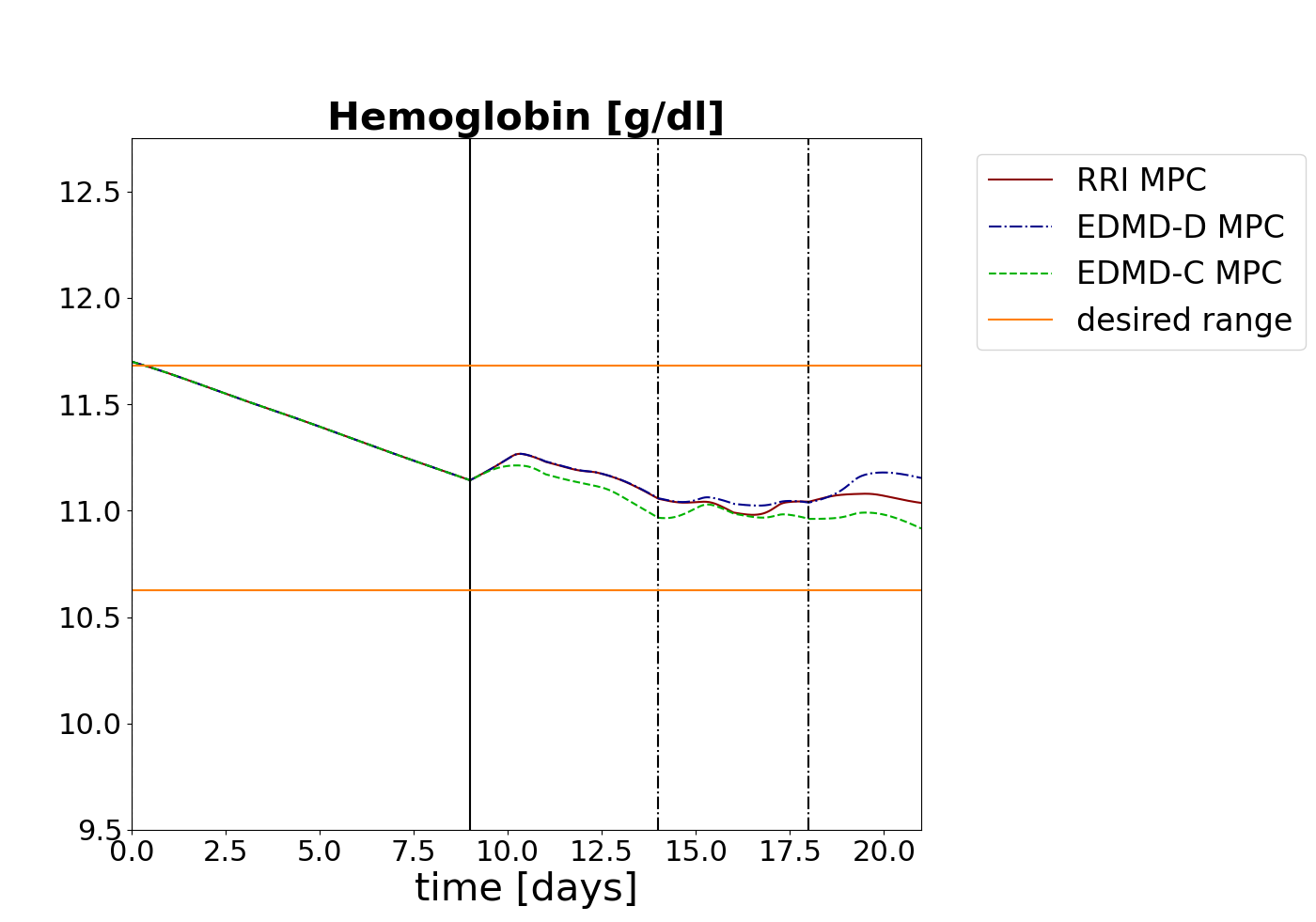}
	\includegraphics[width=0.49\textwidth]%
	{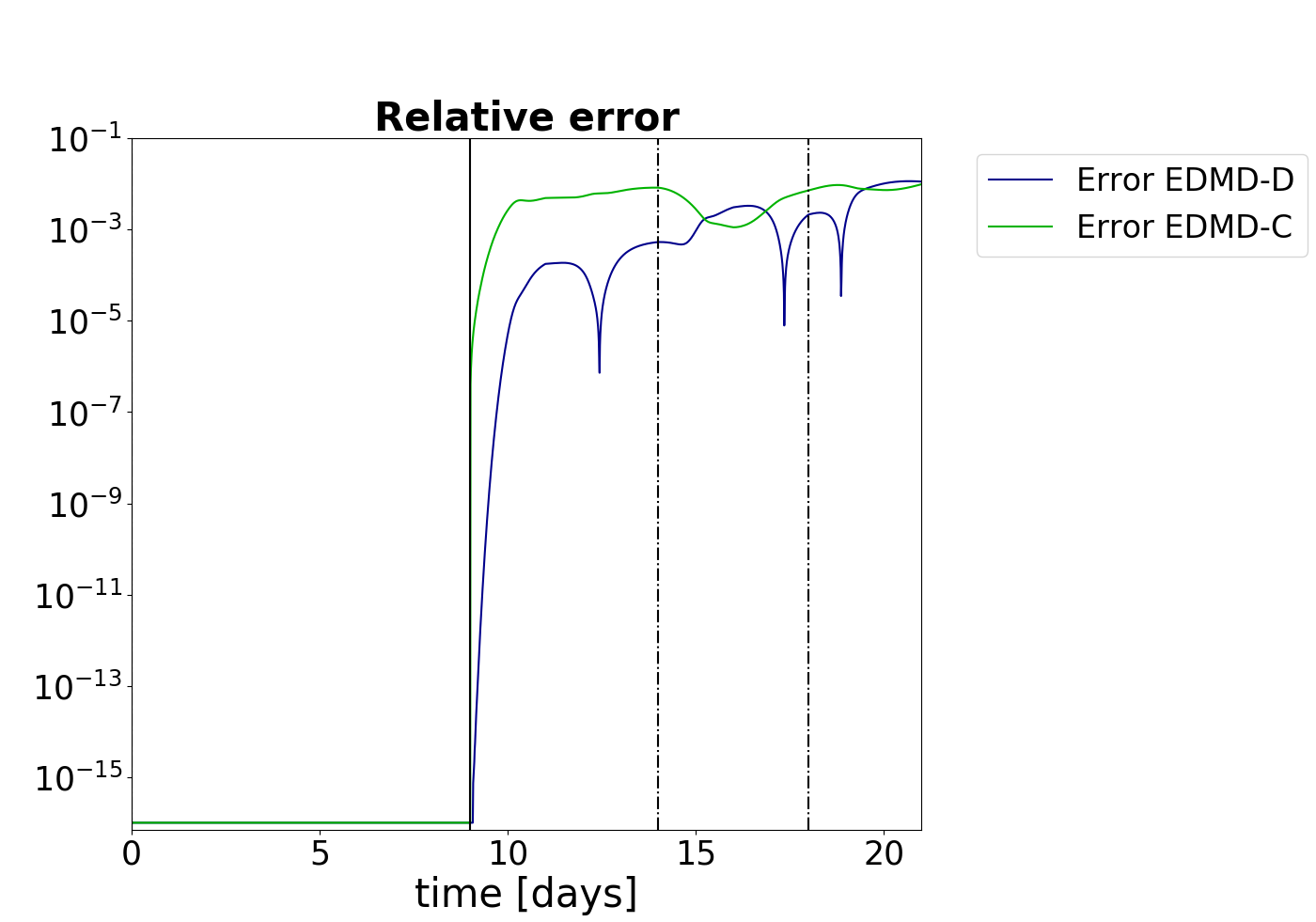} 
	\caption{Patient $2$ - 21 days. {\bf Left:} MPC solutions.
		 {\bf Right:} Relative error. 
		(Black line: EDMD-D update, dashed black line: EDMD-C update). 	\label{Mechelli_MS23_fig:Patient2}}
	\vspace{-1em}
\end{figure}
\begin{figure}
	\vspace{-0.5em}
	\sidecaption
	\includegraphics[width=0.49\textwidth]%
	{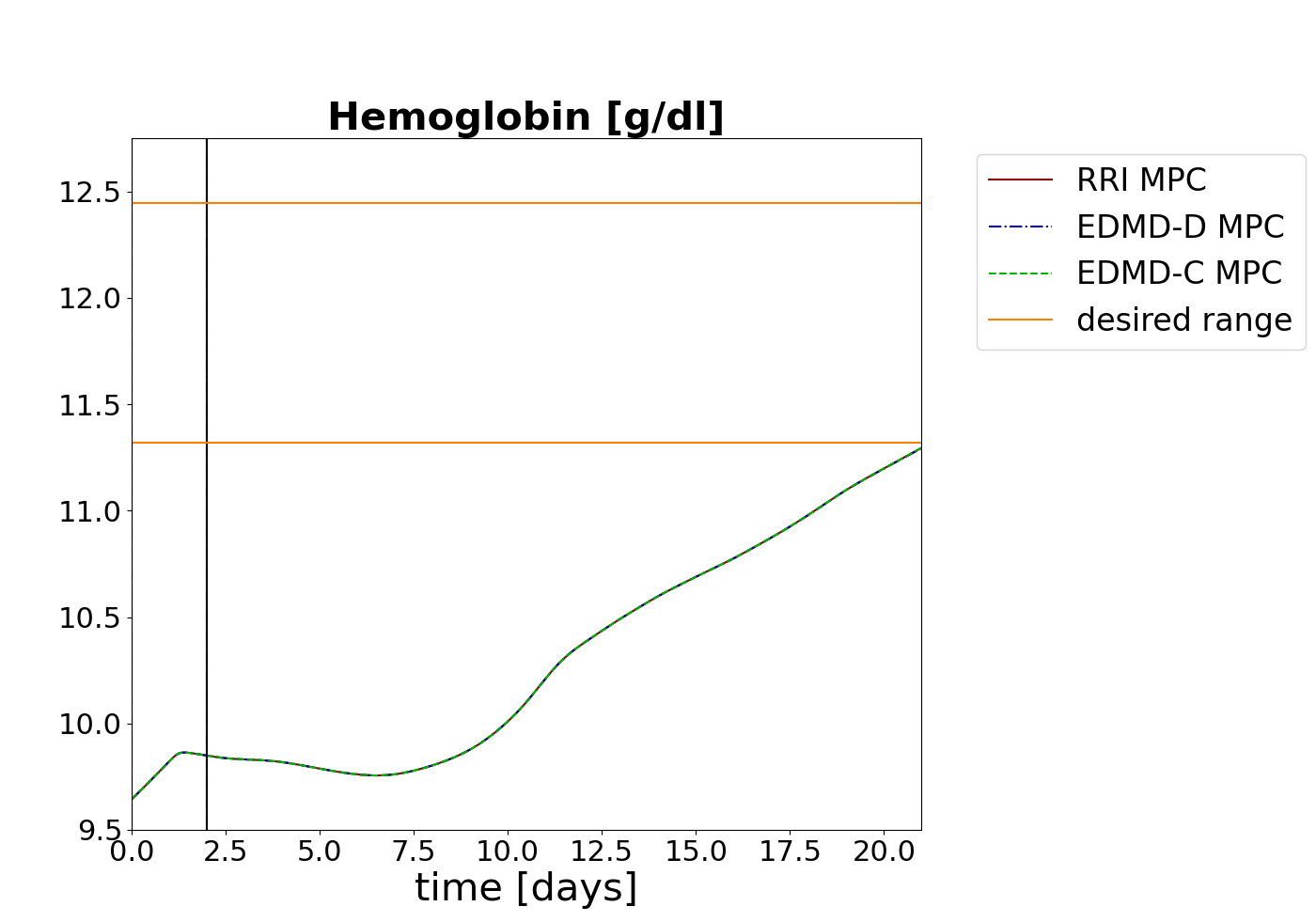}
	\includegraphics[width=0.49\textwidth]%
	{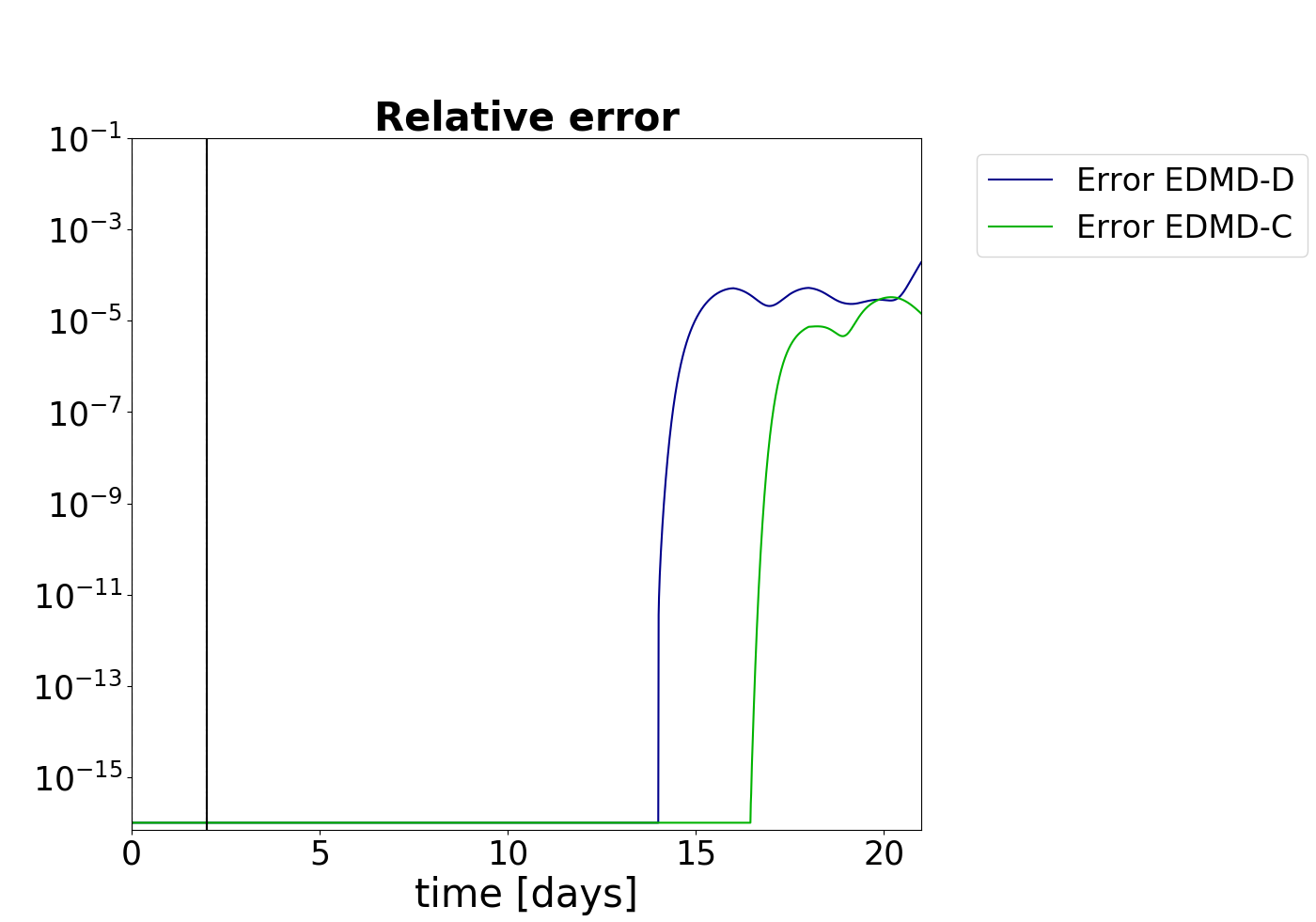}
	\caption{Patient $3$ - 21 days. {\bf Left:} MPC solutions. {\bf Right:} Relative error.  
		(Black line: EDMD-D update, dashed black line: EDMD-C update). \label{Mechelli_MS23_fig:Patient3}}
	\vspace{-1em}
\end{figure}
\begin{figure}
	\vspace{-0.5em}
	\sidecaption
	\includegraphics[width=0.49\textwidth]%
	{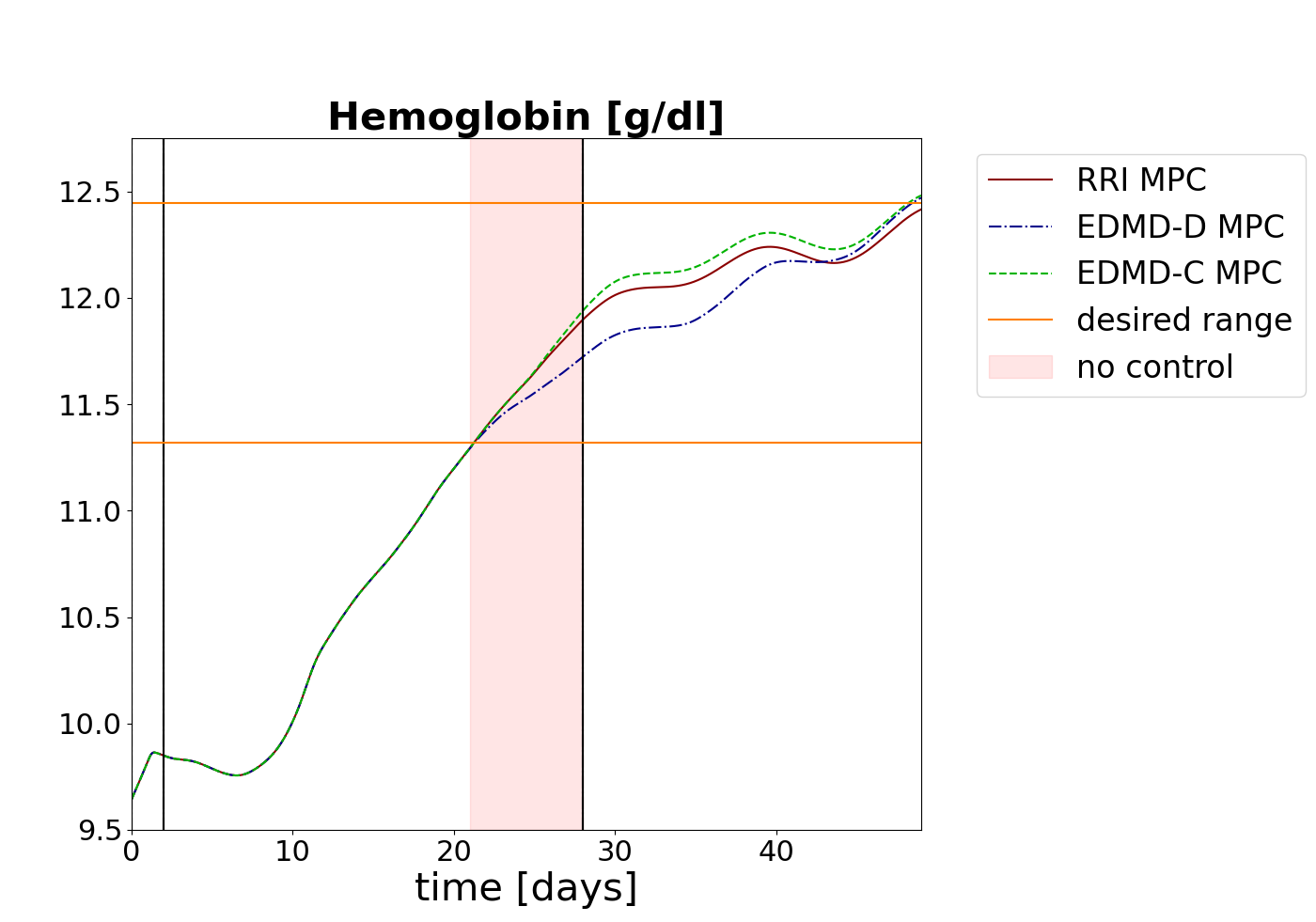}
	\includegraphics[width=0.49\textwidth]%
	{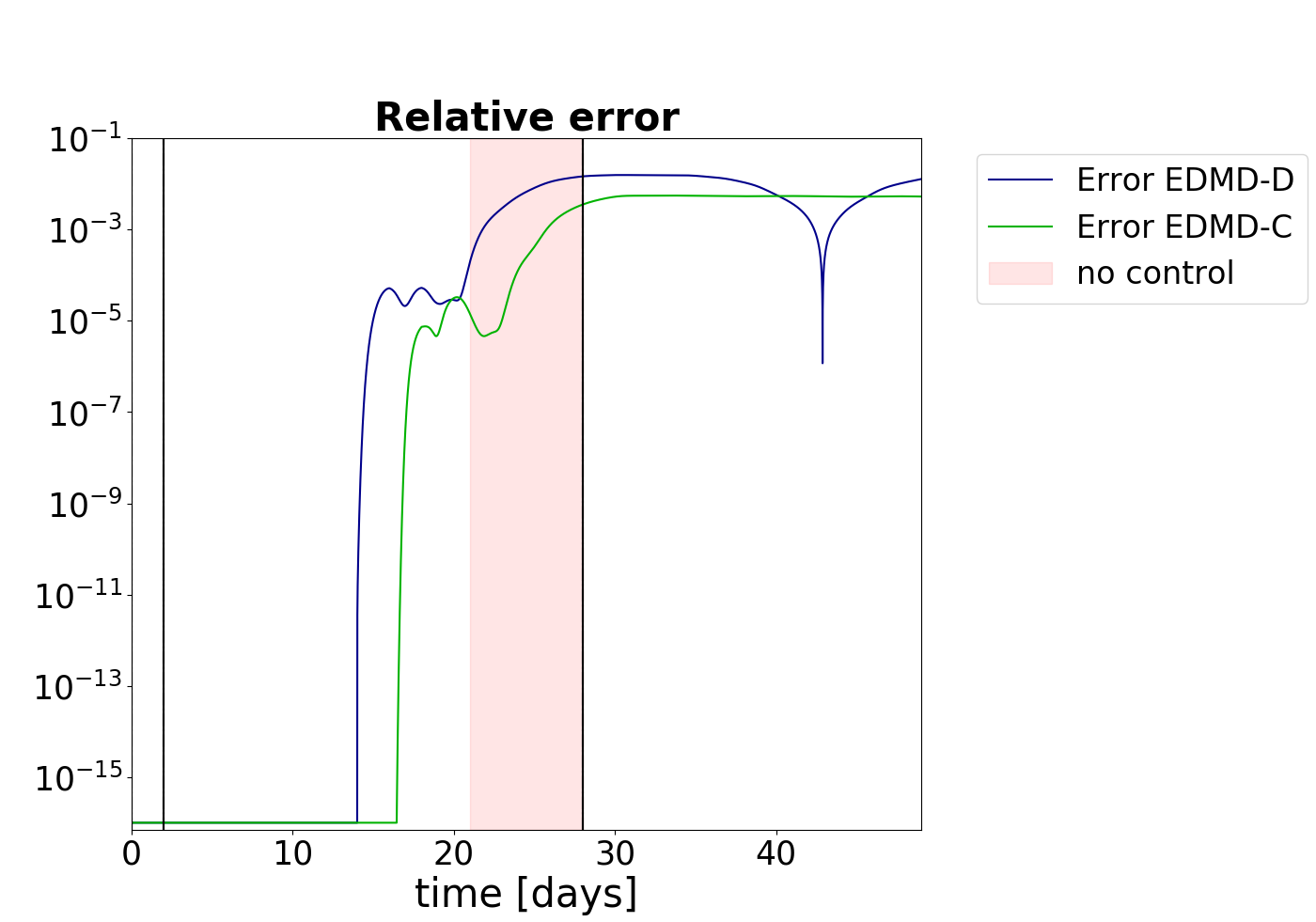}
	\caption{Patient $3$ - 49 days. {\bf Left:} MPC solutions. {\bf Right:} Relative error.  
		(Black line: EDMD-D update, dashed black line: EDMD-C update, Red area: No injection possible). \label{Mechelli_MS23_fig:skipOneWeek}}
	\vspace{-0.5em}
\end{figure}
In Figure~\ref{Mechelli_MS23_fig:Patient2}-left panel, we plot the optimal hemoglobin level computed using the method from \cite{Mechelli_MS23_RFVKK19} (RRI MPC) and with the two EDMD-MPC approaches proposed in Section~\ref{Mechelli_MS23_sec:2}. As one can see, for the first part of the time horizon, the three methods compute exactly the same optimal solution. These corresponds to a doses $u=0$, which is one of the constraint imposed on the doses. As soon as the control starts impacting the system, we can notice some differences arising between the RRI MPC and our proposed method, in particular for the EDMD-C approach. At each time step, such a difference remains approximately smaller than the $1\%$, as it can be seen from Figure~\ref{Mechelli_MS23_fig:Patient2}-right panel, where the relative error between RRI MPC and our method is reported. We observe similar results for Patient 3 as for Patient 2, see Figure~\ref{Mechelli_MS23_fig:Patient3} for further details. It appears that reconstructing the optimal hemoglobin level works extremely well (up to machine precision) for the first half of the time horizon and become worse in the second half (cf. Figure~\ref{Mechelli_MS23_fig:Patient3}). In this case, note that the EDMD updates are not triggered, even though the approximation worsens as time passes. This suggest that the chosen update strategy is too simple and its improvement will be object of future work. In the next test we consider a longer total time period of 7 weeks, i.e. $T=49$ days, for Patient 3. 
\begin{table}[t]
	\vspace{-0.5em}
	\caption{Computational time, speed-up w.r.t. \cite{Mechelli_MS23_RFVKK19} and relative error of the EDMD-MPC algorithm. \label{Mechelli_MS23_tab:1} }
	\begin{tabular}{p{1.5cm}p{2cm}p{1.1cm}p{3.3cm}p{1.6cm}p{1.7cm}}
		\hline\noalign{\smallskip} 
		Patient & Method & $T$ & Computational time \newline (including updates) & Speed-up & Relative error \\
		\noalign{\smallskip}\svhline\noalign{\smallskip}
		2 & EDMD-D & $\text{ }21 \text{ d }$ & $11.5 \text{ s}$ & $\phantom{0}9.5$ & $9.1 \times 10^{-6\phantom{0}}$ \\
		2 & EDMD-C & $\text{ }21 \text{ d }$ & $11.4 \text{ s}$ & $\phantom{0}9.6$  & $1.9\times 10^{-5\phantom{0}}$ \\
		3 & EDMD-D &$\text{ }21 \text{ d }$ & $\phantom{0}5.3 \text{ s}$ & $10.3$ & $9.0 \times 10^{-10}$ \\
		3 & EDMD-C & $\text{ }21 \text{ d }$& $\phantom{0}6.0 \text{ s}$  &$ \phantom{0}9.1$  & $1.0 \times10^{-10}$ \\
		3 & EDMD-D & $\text{ }49 \text{ d }$ & $17.5 \text{ s}$ &  $\phantom{0}7.1$ & $7.4 \times 10^{-5\phantom{0}}$ \\
		3 & EDMD-C & $\text{ }49 \text{ d }$ & $\phantom{0}9.9 \text{ s}$  & $12.6$  & $1.4 \times10^{-5\phantom{0}}$ \\
		\noalign{\smallskip}\hline\noalign{\smallskip}  
	\end{tabular}
\vspace{-0.5em}
\end{table}
This test simulates that Patient 3 skips one week of injections for some reason. In Figure \ref{Mechelli_MS23_fig:skipOneWeek} the week without injections is represented by the red area. For the MPC framework we do not require this additional information. We just assume that the patient will get three injections every week, until he skips the injection appointments. The test demonstrates the ability of the MPC algorithm to react with a quick feedback response. In Figure \ref{Mechelli_MS23_fig:skipOneWeek} we see that our MPC methods based on EDMD are reacting as well (additionally updating the EDMD approximation) and their response is close to the non-linear MPC from \cite{Mechelli_MS23_RFVKK19}. For all the numerical tests, we report the space-time relative error in reconstructing the RRI MPC hemoglobin level and the required computational times for the proposed EDMD-MPC schemes in Table~\ref{Mechelli_MS23_tab:1}. Note that our methods are almost one order of magnitude faster and reconstruct the solution with a reliable approximation error. In conclusion, the proposed EDMD-MPC method replicates the results obtained in \cite{Mechelli_MS23_RFVKK19} with reasonable error and a small factor of speed-up. This factor remains constant when the treatment horizon increases (cf. Table~\ref{Mechelli_MS23_tab:1}). The EDMD-MPC algorithm can be then a valid approach for hemodialysis treatments, although the estimates of the EDMD error and the resulting update strategy during the MPC iterations need to be improved. This will be the focus of a future work.
\begin{acknowledgement}
\vspace{-1.0em}
L. Mechelli and J. Rohleff gratefully acknowledge support of the Independent Research Grant \textit{Efficient model order reduction for model predictive control of non-linear input-output dynamical systems} of the Zukunftskolleg at the University of Konstanz. The authors thank Christian Himpe, Max Planck Institute Magdeburg, for the fruitful discussions and remarks.
\vspace{-2.5em}
\end{acknowledgement}
%
%
%

\begin{thebibliography}{99.}%
%
%
\bibitem{Mechelli_MS23_Fue12}
Fuertinger, D.H.: A model of erythropoiesis. PhD thesis, Karl-Franzens University Graz (2012)

\bibitem{Mechelli_MS23_FKTLK13}
Fuertinger, D.H., Kappel, F., Thijssen, S., Levin, N.W., Kotanko, P.: A model of erythropoiesis in adults with sufficient iron availability. J. Math. Biol. \textbf{66}(6), 1209--1240 (2013)

\bibitem{Mechelli_MS23_GP17}
Gr\"une, L., Pannek, J.: Nonlinear Model Predictive Control:Theory and Algorithms. 2nd Edition. Springer, London (2016)

\bibitem{Mechelli_MS23_HPUU09}
Hinze, M., Pinnau, R., Ulbrich, M., Ulbrich, S.:  Optimization with PDE Constraints. Springer-Verlag, Berlin (2009)

\bibitem{Mechelli_MS23_KM2018}
Korda, M., Mezi\'c I.: Linear predictors for nonlinear dynamical systems: Koopman oprator meets model predictive control. Automatica \textbf{93}, 149--160 (2018)

\bibitem{Mechelli_MS23_KBBP16}
Kutz, J.N., Brunton, S.L., Brunton, B.W., Proctor J.L.: Dynamic Mode Decomposition: Data-Driven Modeling of Complex Systems. SIAM, Philadelphia (2016)

\bibitem{Mechelli_MS23_RFVKK19}
Rogg, S., Fuertinger, D.H., Volkwein, S., Kappel, F., Kotanko, P.: Optimal EPO dosing in hemodialysis patients using a non-linear model predictive control approach. J. Math. Biol. \textbf{79}, 2281--2313 (2019)

\bibitem{Mechelli_MS23_WKR15}
Williams, M.O., Kevrekidis, I.G., Rowley, C.W.: A Data-Driven Approximation of the Koopman Operator: Extending {D}ynamic {M}ode {D}ecomposition. J. Nonlinear Sci. \textbf{25}, 1307--1346 (2015)

\end{thebibliography}
%

\end{document}